\documentclass[a4paper,twoside,
11pt]{amsart}
\usepackage{amssymb,amsfonts,amsmath}
\begin{document}
\newtheorem{defin}{~~~~Definition}
\newtheorem{prop}{~~~~Proposition}[section]
\newtheorem{remark}{~~~~Remark}[section]
\newtheorem{cor}{~~~~Corollary}
\newtheorem{theor}{~~~~Theorem}
\newtheorem{lemma}{~~~~Lemma}[section]
\newtheorem{ass}{~~~~Assumption}
\newtheorem{con}{~~~~Conjecture}
\newtheorem{concl}{~~~~Conclusion}
\numberwithin{equation}{section}
\newcommand {\trans} {^{\,\mid\!\!\!\cap}}
\newcommand{\vf}{\varphi}
\newcommand{\mJ}{\mathcal J}
\newcommand{\e}{\varepsilon}
\title{ A Canonical Frame for Nonholonomic
Rank Two Distributions of Maximal Class}
\author
{Boris Doubrov
\address{
the Faculty of Applied Mathematics, Belorussian State
University, F. Skaryny Ave. 4, Minsk, Belarus 220050;
 E-mail: doubrov@islc.org}\and Igor Zelenko
\address
{
S.I.S.S.A., Via Beirut 2-4, 34014, Trieste,
Italy;
E-mail:
 zelenko@sissa.it}} \maketitle \markboth{Boris
Doubrov and Igor Zelenko} {A Canonical Frame for Rank 2
Distributions}

\begin{abstract}
In 1910 E. Cartan constructed the canonical frame and found
the most symmetric case for maximally nonholonomic rank $2$
distributions in $\mathbb R^5$.
We
solve the analogous problems for rank $2$ distributions in
${\mathbb R}^n$ for arbitrary $n>5$.
Our method is a kind of symplectification of the problem
and it is completely different from the Cartan method of
equivalence.
\end{abstract}

\section{Introduction}

 A
rank $l$ vector distribution $D$ on an $n$-dimensional
manifold $M$ or an $(l,n)$-distribution (where $l<n$) is a
subbundle of the tangent bundle $TM$ with $l$-dimensional
fibers. The group of germs of diffeomorphisms of $M$ acts
naturally on the set of germs of $(l,n)$-distributions and
defines the equivalence relation there.
The question is \emph{when two germs of distributions are 
equivalent?} Distributions are naturally associated with 
Pfaffian systems and with control systems linear in the 
control. So the problem of equivalence of distributions can 
be reformulated as the problem of equivalence of the 
corresponding Pfaffian systems and the state-feedback 
equivalence of the corresponding control systems. The 
obvious (but very rough in the most cases) discrete 
invariant of a distribution $D$ at $q$ is so-called \emph{ 
the small growth vectors} at $q$. It is the tuple 
$\{\dim D^j(q)\}_{j\in{\mathbb N}}$, where $D^j$ is the 
$j$-th power of the distribution $D$, i.e., 
$D^j=D^{j-1}+[D,D^{j-1}]$, $D^1=D$. 
A simple estimation shows that at least $l(n-l)-n$ 
functions of $n$ variables are required to describe generic 
germs of $(l,n)$-distribution, up to the equivalence (see 
\cite{versh1} and \cite {zhit0} for precise statements). 
 There are 
only three cases, where $l(n-l)-n$ is not positive: $l=1$ 
(line distributions), $l=n-1$, and $(l,n)=(2,4)$. Moreover, 
it is well known that in these cases generic germs of 
distributions are equivalent.
for $l=1$ it is just the classical theorem about the 
rectification of vector fields without stationary points, 
for $l=n-1$ all generic germs are equivalent to Darboux's 
model, while for $(l,n)=(2,4)$ they are equivalent to 
Engel's model (see, for example, \cite{bryantbook}). In all 
other cases generic $(l,n)$-distributions have functional 
invariants. 

In the present paper we restrict ourselves to the case of 
rank 2 distributions. The model examples of such 
distributions come from so-called underdetermined ODE's of 
the type 
$$z^{(r)}(x)=F\bigl(x,y(x),\ldots, y^{(s)}(x),z(x),\ldots, 
z^{(r-1)}(x)\bigr), \quad r+s=n-2,$$ 
for two functions $y(x)$ and $z(x)$. Setting $p_i=y^{(i)}$, 
$0\leq i\leq s$, and $q_j=z^{(j)}$, $0\leq j\leq r-1$, with 
each such equation one can associate the rank $2$ 
distribution 
in ${\mathbb R}^n$ with coordinates $(x,p_0,\ldots, 
p_s,q_0,\ldots, q_{r-1})$ given by the intersection of the 
annihilators of the following $n-2$ one-forms: 
\begin{equation*}
\begin{aligned} 
~&dp_i-p_{i+1} dx , \,\,0\leq i\leq s-1, \quad 
 dq_j-q_{j+1} dx,\,\, 0\leq j\leq r-2,\\
 ~&dq_{r-1}-F(x,p_0,\ldots, p_s,q_0,\ldots, q_{r-1})dx.
\end{aligned}
\end{equation*}

%
%
%
%




For $n=3$ and $4$ all generic germs of rank 2 distribution 
are equivalent to the distribution, associated with the 
underdetermined ODE $z'(x)=y(x)$ (Darboux and Engel models 
respectively). 
The case 
$n=5$ (the smallest dimension, when functional parameters 
appear) 
was treated by E. Cartan in \cite{cartan} with his 
reduction-prolongation procedure. First, for any 
$(2,5)$-distribution with the small growth vector $(2,3,5)$ 
he constructed the canonical coframe in some 14-dimensional 
manifold, which implied that the group of symmetries of 
such distributions is at most $14$-dimensional. Second, 
he showed that any $(2,5)$-distribution with 
$14$-dimensional group of symmetries is locally equivalent 
to the distribution, associated with the underdetermined 
ODE $z'(x)=\bigl(y''(x)\bigr)^2$, and its group of 
symmetries 
is isomorphic to the real 
split form of the exceptional Lie group $G_2$. Historically 
it was the first natural appearance of this group. 

After the work of Cartan the open question was \textbf {to 
construct the canonical frame and to find the most 
symmetric cases for $(2,n)$-distributions with $n>5$}.
The Cartan equivalence method was systematized and 
generalized by N. Tanaka and T. Morimoto (see \cite{tan2, 
mori}). Their theory is heavily based on the notion of 
so-called \emph{ symbol algebra} of the distribution at a 
point, which is a special \emph{ graded nilpotent Lie 
algebra}, naturally associated with the distribution at a 
point: the symbol algebras have to be isomorphic at 
different points and all constructions strongly depend on 
the type of the symbol. Note that already in the case of 
$(2,6)$-distributions with maximal possible small growth 
vector $(2,3,5,6)$ three different symbol algebras are 
possible, while for $n=9$ the set of all possible symbol 
algebras depends on continuous parameters, which implies in 
particular that generic distributions do not have a 
constant symbol. 

In the present paper we give an answer to the question, 
underlined in the previous paragraph, for rank $2$ 
distributions from some generic class. Our constructions 
are based on a completely different, variational approach, 
developed in \cite{jac1} and \cite{zelvar}. 
Roughly speaking, we make a kind of symplectification of 
the problem by lifting the distribution to the cotangent 
bundle $T^*M$ of the 
manifold $M$. 

\section {The class of rank 2 distribution}
Assume that $\dim D^2(q)=3$ and 
$\dim D^3(q)>3$ for any $q\in M$. Denote by
$(D^l)^{\perp}\subset T^*M$ the annihilator of the $l$th 
power $D^l$, namely
$$(D^l)^{\perp}= \{(q,p)\in T^*M:\,\, p\cdot v=0\,\,\forall 
v\in D^l(q)\}.$$ 

First, we distinguish a characteristic $1$-foliation on the 
codimension $3$ submanifold 
$(D^2)^\perp\backslash(D^3)^\perp$ of $T^*M$. For this let 
$\pi\colon T^*M\mapsto M$ be the canonical projection. For 
any $\lambda\in T^*M$, $\lambda=(p,q)$, $q\in M$, $p\in 
T_q^*M$, let $\mathfrak{s}(\lambda)(\cdot)=p(\pi_*\cdot)$ 
be the canonical Liouville form and $\sigma=d\mathfrak {s}$ 
be the standard symplectic structure on $T^*M$. Since the 
submanifold $(D^2)^\perp$ has odd codimension in $T^*M$, 
the kernels of the restriction 
$\sigma|_{(D^2)^\perp}$ of $\sigma$ on $(D^2)^\perp$ are 
not trivial. Moreover for 
the points of $(D^2)^\perp\backslash (D^3)^\perp$ these 
kernels are one-dimensional. 
They form the \emph{characteristic line distribution} in 
$(D^2)^\perp\backslash(D^3)^\perp$, which will be denoted 
by ${\mathcal C}$. The line distribution ${\mathcal C}$ 
defines a \emph{characteristic 1-foliation} of 
$(D^2)^\perp\backslash(D^3)^\perp$. The leaves of this 
foliation are called the \emph{characteristic curves}.
 In Control Theory these characteristic 
curves are also called \emph{regular abnormal extremals of 
$D$}. 

In the sequel given two submanifold $S_1$ and $S_2$ of the 
tangent bundle of some manifold $W$ such that 
$S_i(w)=S_i\cap T_w W$, $i=1,2$, are linear subspaces of 
$T_wW$ (not necessary of the same dimensions for different
$w$) we will denote by $[S_1, S_2]$ the subset $\{[S_1, 
S_2](w)\}_{w\in W}$ of $TW$ such that $$[S_1, S_2](w)={\rm 
span}\{[Z_1, Z_2](w):Z_i\,\,{\rm are}\,\, {\rm vector}\,\, 
{\rm fields}\,\, {\rm tangent}\,\,{\rm to}\,\, S_i, 
i=1,2\}.$$ It is easy to show that with such definition 
$S_i\subset [S_1, S_2]$ , $i=1,2$. 


Now, following \cite{zelvar}, let  
\begin{equation}
\label{prejac} 
{\mathcal J}(\lambda)= 
\bigl(T_\lambda 
(T^*_{\pi(\lambda)}M)+ 
\ker\sigma|_{D^\perp}(\lambda)\bigr)\cap T_\lambda 
(D^2)^\perp= 
\{v\in T_{\lambda}(D^2)^\perp:\,\pi_*\,v\in 
D(\pi\bigl(\lambda)\bigr)\}.
\end{equation} 
(here $T_\lambda (T^*_{\pi(\lambda)}M)$ is tangent to the 
fiber $T^*_{\pi(\lambda)}M$ at the point $\lambda$). 
Note that $\dim {\mathcal J}(\lambda)=n-1$. 
Actually, ${\mathcal J}$ is the pull-back of the 
distribution $D$ on $(D^2)^\perp\backslash(D^3)^\perp$ by 
the canonical projection $\pi$. Another important property of the distribution 
$\mJ$, which follows from its definition, is 
$\sigma|_\mJ=0$.
Define a sequence of subspaces ${\mathcal 
J}^{(i)}(\lambda)$, $\lambda\in (D^2)^\perp\backslash 
(D^3)^\perp$, by the following recursive formulas: 
\begin{equation}
\label{Ji1} {\mathcal J}^{(i)}=
[{\mathcal C},{\mathcal J}^{(i-1)}] \quad {\mathcal J}^{(0)}={\mathcal J}.
\end{equation}
 By \cite{zelvar}( Proposition 
3.1 and formula (3.9) 
there), 
\begin{equation}
\label{equivrk} 
{\rm dim}\,{\mathcal J}^{(1)}(\lambda)-{\rm dim}\, 
{\mathcal J}(\lambda) =1 
\end{equation}
Actually, \begin{equation} \label{J1} {\mathcal 
J}^{(1)}(\lambda)=\{v\in 
T_{\lambda}(D^2)^\perp:\,\pi_*\,v\in 
D^2(\pi\bigl(\lambda)\bigr)\}. 
\end{equation}

Relation (\ref{equivrk}) together with the definition 
(\ref{Ji1}) implies that 
\begin{equation}
\label{difdimJ} 
{\rm dim}\, {\mathcal J}^{(i)}(\lambda)- {\rm dim}\, 
{\mathcal J}^{(i-1)}(\lambda)\leq 1,\quad  i\in {\mathbb 
N}.
\end{equation}

\begin{lemma}
\label{estlemma} The following inequality holds
\begin{equation}
\label{est} 
 {\rm dim}\, {\mathcal 
J}^{(i)}(\lambda)\leq 2n-4. 
\end{equation}
\end{lemma}

{\bf Proof.} By definition the line distribution ${\mathcal 
C}$ forms the characteristic of the corank $1$ distribution 
on $(D^2)^\perp$, given by the Pfaffian equation $\mathfrak 
{s}|_{(D^2)^\perp}=0$. Since by construction ${\mathcal 
J}\subset\{{\mathfrak s}|_{(D^2)^\perp}=0\}$, one has 
\begin{equation}
\label{inLiuv} {\mathcal J}^{(i)}\subset\{{\mathfrak 
s}|_{(D^2)^\perp}=0\}\quad i\in\mathbb N. 
\end{equation}
Our lemma follows from the fact that the distribution 
$\{{\mathfrak s}|_{(D^2)^\perp}=0\}$ has rank $2n-4$. 
$\Box$ 
\medskip

Further for any point $q\in M$ 
denote by $(D^l)^\perp(q) =(D^l)^\perp\cap T_q^*M$. Let us 
define the following two integer-valued functions: 
\begin{equation*}
\nu(\lambda)=\min\{i\in {\mathbb N}: {\mathcal 
J}^{(i+1)}(\lambda)={\mathcal J}^{(i)}(\lambda)\}, 
\quad
m(q)=\max\{\nu(\lambda):\lambda\in 
(D^2)^\perp(q)\backslash(D^3)^\perp(q)\}. 
\end{equation*}
\begin{defin}
\label{classdef} The number
$m(q)$ is 
called the class of distribution $D$ at the point $q$. 
\end{defin}

By (\ref{equivrk}), (\ref{difdimJ}), and the previous lemma 
$1\leq m(q)\leq n-3$. It is easy to show that the 
integer-valued functions $\nu(\cdot)$ and $m(\cdot)$ are 
lower semicontinuous. Hence they are locally constant on 
the open and dense subset of $(D^2)^\perp\backslash 
(D^3)^\perp$ and $M$ correspondingly. We say that the point 
$q\in M$ is \emph{the regular point} of the distribution 
$D$, if the function $m(\cdot)$ is constant in some 
neighborhood $U$ of $q$, i.e. the distribution $D$ has 
constant class on $U$. Obviously, all points, where the 
distribution $D$ has maximal class $n-3$, are regular. 
Moreover, 

\begin{prop}
\label{nonhol} 
 Germs of $(2,n)$-distributions of the maximal class $n-3$  
 are generic. 
\end{prop}
This Proposition follows directly from Proposition 3.4 of 
\cite{zelvar}. In the present paper we treat the germs of 
$(2,n)$ distributions of the maximal class $n-3$. In the 
cases $n=5$ and $n=6$ a rank 2 distribution has maximal 
class if and only if it has maximal possible small growth 
vector, namely, $(2,3,5)$ in the case $n=5$ and $(2,3,5,6)$ 
in the case $n=6$ (see Propositions 3.5 and 3.6 of 
\cite{zelvar} respectively). 



\section{The canonical projective structure on characteristic curves.} 
From now on $D$ is a
$(2,n)$-distribution of maximal constant class $m=n-3$. Let 
\begin{equation*} 
\label{Jconstw2} 
{\mathcal R}_D=\{\lambda\in (D^2)^\perp \backslash 
(D^3)^\perp: \nu(\lambda)=n-3 \}, \quad  {\mathcal 
R}_D(q)={\mathcal R}_D\cap T_q^*M. 
\end{equation*}
Note that the set ${\mathcal R}_D(q)$ is a nonempty open 
set in Zariski topology on the linear space 
$(D^2)^\perp(q)$ (see again~\cite[Proposition~3.4]{zelvar}).
The crucial observation is that \emph{any segment of 
a characteristic curve $\gamma$ of $D$, belonging to 
${\mathcal R}_D$, can be endowed with a canonical 
projective structure} (for more detailed description than 
below see \cite{agrzel},\cite{jac1}, and \cite{zelvar}). By 
a projective structure on a curve we mean a set of 
parameterizations such that the transition function from 
one such parameterization to another is a M\"{o}bious 
transformation. To construct this canonical projective 
structure on $\gamma$ first we associate with $\gamma$ a 
special curve in a Grassmannian $G_m(W)$ of $m$-dimensional 
subspaces of a $2m$-dimensional linear space $W$, 
the Jacobi curve, in the following way: Let 
$O_\gamma$ be a neighborhood of $\gamma$ in $(D^2)^\perp$ 
such that the factor 
$N=O_\gamma /(\text {\emph{the characteristic 
one-foliation}})$ 
 is a 
well-defined smooth manifold. Its dimension is equal to 
$2(n-2)$. 
Let $\phi \colon O_\gamma\to N$ be the canonical projection 
on the factor. 
Define the mapping  
$J_\gamma\colon\gamma\mapsto G_{n-2}(T_\gamma N)$ by 
$J_\gamma(\lambda) 
{=} \phi_*\bigl({\mathcal J}(\lambda) \bigr)$ for all 
$\lambda\in\gamma$, where ${\mathcal J}(\lambda)$ is as 
above. Actually, the symplectic form $\sigma$ of $T^*M$ 
induces naturally the symplectic form $\bar \sigma$ on 
$T_\gamma N$ and $J_\gamma(\lambda)$ for all
$\lambda\in\gamma$ are Lagrangian subspace of $T_\gamma N$. 
Besides, if $e$ is the Euler field (i.e., the infinitesimal 
generator of homothethies on the fibers of $T^*M$), then 
the vector $\bar e= \phi_*e(\lambda)$ is the same for any 
$\lambda\in \gamma$ and lies in $J_\gamma(\lambda)$. 
Therefore the curve $\lambda\mapsto\widetilde 
J_\gamma(\lambda)=J_\gamma(\lambda)/\{\mathbb R \bar e\}$, 
$\lambda\in \gamma$, is a curve in $G_{m}(W)$, where 
$W=\{v\in T_\gamma N: \sigma(v,\bar e)=0\}/\{\mathbb R \bar 
e\}$
. The curve $\widetilde J_\gamma$ is called \emph{ the 
Jacobi curve of $\gamma$}. 

Second, we construct the canonical projective structure on 
$\widetilde J_\gamma$ (and therefore on $\gamma$ itself), 
using the notion of the generalized cross-ratio of $4$ 
points in $G_{m}(W)$. Namely, let 
$\{\Lambda_i\}_{i=1}^4$ be any $4$ points of $G_{m}(W)$. 
For simplicity suppose that $\Lambda_i\cap\Lambda_j=0$ for 
$i\neq j$. Assume that in some coordinates $W\cong \mathbb 
R^m\times\mathbb R^m$ and $\Lambda_i=\{(x, S_ix):x\in 
\mathbb R^m\}$ for some $m\times m$-matrix $S_i$. Then the 
conjugacy class of the following matrix 
$$(S_1-S_4)^{-1}(S_4-S_3) (S_3-S_2)^{-1}(S_2-S_1)$$ does 
not depend on the choice of the coordinates in $W$. This 
conjugacy class is called \emph{the cross-ratio of the 
tuple $\{\Lambda_i\}_{i=1}^4$} and it is denoted by 
$[\Lambda_1,\Lambda_2,\Lambda_3,\Lambda_4]$. 

Now take some parametrization $\varphi\colon\gamma\mapsto 
{\mathbb R}$ of $\gamma$ and let 
$\Lambda_\varphi(t)=\widetilde 
J_\gamma\bigl(\varphi^{-1}(t)\bigr)$. Assume that in some 
coordinates on $W$ we have 
$\Lambda_\varphi(t)=\{(x,S_tx):x\in\mathbb R^m\}$. 
The following fact follows from \cite[Proposition 
2.1]{zelvar} : For all parameters $t_1$ the functions 
$t\rightarrow det(S_t-S_{t_1})$ have zero of the the same 
order $k=m^2$ at $t=t_1$. Consider the following function 
\begin{equation}
\label{fung}{\mathcal G}_{\Lambda_\varphi} 
(t_1,t_2,t_3,t_4)=ln 
\left(\det[\Lambda_\varphi(t_1),\Lambda_\varphi(t_2), 
\Lambda_\varphi(t_3),\Lambda_\varphi(t_4)]\bigl([t_1,t_1,t_2,t_3 
]\bigr)^{-k} \right),
\end{equation}
 where 
$[t_1,t_2,t_3,t_4]=\frac{(t_1-t_2)
(t_2-t_3)}{(t_3-t_4)(t_1-t_4)}$ is the usual cross-ratio of 
$4$ numbers $\{t_i\}_{i=1}^4$. Then, by above, it is not
hard to see that ${\mathcal G}_{\Lambda_\varphi} 
(t_1,t_2,t_3,t_4)$ is smooth at diagonal points $(t,t,t,t)$ 
and the Taylor expansions up to the order two of it at 
these points have the form 
\begin{equation}
\label{rho} {\mathcal G}_{\Lambda_\varphi} 
(t_0,t_1,t_2,t_3)=\rho_{\Lambda_\varphi}(t) 
(\xi_1-\xi_3)(\xi_2-\xi_4)+\ldots,\quad \xi_i=t_i-t.
\end{equation}
Now let $\psi:\mathbb R\mapsto \mathbb R$ be a smooth 
monotonic function. 
Then by (\ref{fung})  
\begin{equation}
\label{repG} {\mathcal 
G}_{\Lambda_\vf}(t_1,t_2,t_3,t_4)={\mathcal 
G}_{\Lambda_{\psi\circ\vf}}\bigl(\psi(t_1),\psi(t_2),\psi(t_3),\psi(t_4)\bigr)+k\, 
\ln\left(\frac {[\psi(t_1),\psi(t_2),\psi(t_3),\psi(t_4)]} 
{[t_1,t_2,t_3,t_4]}\right), 
\end{equation}
By direct computation it can be shown that the function 
$(t_0,t_1,t_2,t_3)\mapsto \ln\left(\frac 
{[\psi(t_1),\psi(t_2),\psi(t_3),\psi(t_4)]} 
{[t_1,t_2,t_3,t_4]}\right)$ has the following Taylor 
expansion up to the order two at the point $(t,t,t,t)$: 
\begin{equation}
\label{lnsch} \ln\left(\frac 
{[\psi(t_1),\psi(t_2),\psi(t_3),\psi(t_4)]} 
{[t_1,t_2,t_3,t_4]}\right)=\cfrac{1}{3}\mathbb 
S\psi(t)(\xi_1-\xi_3)(\xi_2-\xi_4)+\ldots,\quad 
\xi_i=t_i-t, 
\end{equation}
where $\mathbb S\psi$ is Schwarz derivative of $\psi$, 
$\mathbb S\psi= 
\frac{1}{2}\frac{\psi^{(3)}}{\psi'}-\frac{3}{4}\Bigl 
(\frac{\psi''}{\psi'}\Bigr)^2$. Combining (\ref{rho}) and 
(\ref{lnsch}) we get the following reparameterization rule 
for $\rho_{_{\Lambda_\vf}}$: \begin{equation} 
\label{reparam} \rho_{_{\Lambda_\vf}}(t)=
\rho_{_{\Lambda_{\psi\circ\vf}}}(\psi(t))(\psi'(t))^2+\cfrac{k}{3}
\,\mathbb S\psi(t). \end{equation} 
From the last formula and the fact that ${\mathbb
S}\psi\equiv 0$ if and only if the function $\psi$ is 
M\"{o}bius it follows that \emph{the set of all 
parametrizations $\varphi$ of $\gamma$ such that 
$\rho_{\Lambda_\varphi}(t)\equiv 0$ defines the canonical 
projective structure on $\gamma$}.


\section{The main theorem}
Now we are ready to describe the 
manifold, on which the canonical frame for 
$(2,n)$-distribution of maximal class, $n>5$, can be
constructed. Given $\lambda\in {\mathcal R}_D$ denote by 
${\mathfrak P}_\lambda$ the set of all projective 
parameterizations $\varphi\colon\gamma\mapsto\mathbb R$ on 
the characteristic curve $\gamma$ , passing through 
$\lambda$, such that $\varphi(\lambda)=0$. Denote 
$$\Sigma_D=\{(\lambda, \varphi):\lambda\in {\mathcal R}_D, 
\varphi\in {\mathfrak P}_\lambda\}.$$ Actually, $\Sigma_D$ 
is a principal bundle over ${\mathcal R}_D$ with the 
structural group of all M\"{o}bious transformations, 
preserving $0$ and $\dim \, \Sigma_D=2n-1$. 

{\bf Theorem.} 
\emph{For any $(2,n)$-distribution, $n>5$, of maximal
class there exist two canonical frames on the corresponding 
$(2n-1)$-dimensional manifold $\Sigma_D$, obtained one from 
another by a reflection. The group of symmetries of such 
distribution is at most $(2n-1)$-dimensional. Any 
$(2,n)$-distribution of maximal class with 
$(2n-1)$-dimensional group of symmetries is locally
equivalent to the distribution, associated with the 
underdetermined ODE $z'(x)=\bigl(y^{(n-3)}(x)\bigr)^2$. The 
algebra of infinitesimal symmetries of this distribution is 
isomorphic to a semidirect sum of $\mathfrak{gl}(2,\mathbb 
R)$ and $(2n-5)$-dimensional Heisenberg algebra ${\mathfrak 
n}_{2n-5}$.} 

{\bf Sketch of the proof.} Define the following two 
fiber-preserving flows on $\Sigma_D$: 
\begin{equation*}
F_{1,s}(\lambda,\vf)=(\lambda, e^{2s}\varphi),\quad 
F_{2,s}(\lambda,\varphi)=\left(\lambda, 
\cfrac{\varphi}{s\varphi+1}\right),\quad  \lambda\in{\mathcal 
R}_D, \varphi\in {\mathfrak P}_\lambda.
\end{equation*}
 
 Further, let $\delta_s$ be the flow of homotheties on 
 the fibers of $T^*M$: 
$\delta_s(p,q)=(e^sp,q)$, where $q\in M$, 
$p\in T_q^*M$ (actually the Euler field $e$ generates this 
flow). The following flow 
$$F_{0,s}(\lambda,\varphi)=\bigl(\delta_s(\lambda), 
\varphi\circ\delta_s^{-1}\bigr)$$ is well-defined on 
$\Sigma_D$ (here we use that $\delta_s$ preserves the 
characteristic $1$-foliation). 
For any $0\leq i\leq 2$ let $g_i$ be the vector field on 
$\Sigma_D$, generating the flow $F_{i,s}$. 
Besides, the characteristic $1$-foliation on $(D^2)^\perp$ 
can be lifted to the \emph{parameterized} $1$-foliation on
$\Sigma_D$, which gives one more canonical vector field on 
$\Sigma_D$. Indeed, let $u=(\lambda,\varphi)\in \Sigma_D$ 
and $\gamma$ be the characteristic curve, passing through 
$\lambda$ (so, $\vf$ maps $\gamma$ to $\mathbb R$). Then 
the the mapping $$\Gamma_u(t)=\bigl(\vf^{-1}(t), 
\vf(\cdot)-t\bigr)$$ defines the parametrized curve on 
$\Sigma_D$, the lift of $\gamma$ to $\Sigma_D$, and 
$\Gamma_u(0)=u$. The additional canonical vector field $h$ 
on $\Sigma_D$ is defined by 
$h(u)=\frac{d}{dt}\Gamma_u(t)|_{t=0}$. It can be shown 
easily that 
\begin{equation}
\label{gl2} 
[g_1,g_2]=2g_2,\,\, [g_1, h]=-2h,\,\, [g_2,h]=g_1,\,\, 
[g_0, h]=0,\,\,[g_0,g_i]=0  
\end{equation}
Therefore 
the linear span (over ${\mathbb R}$) of the vector fields 
$g_0$, $g_1$, $g_2$, and $h$ is endowed with a structure of the  
Lie algebra isomorphic to 
$\mathfrak{gl}(2,\mathbb{R})$. 

Now we will construct one more canonical, up to the sign, 
vector field on $\Sigma_D$. 
For this let 
\begin{equation}
\label{contr} {\mathcal J}_{(i)}(\lambda)= \{v\in T_\lambda 
\bigl((D^2)^\perp\bigr): \sigma (v, w)=0\,\, \forall w\in 
{\mathcal J}^{(i)}\},\quad V_i(\lambda)=\{\lambda\in 
\mJ_{(i)}:\pi_*(v)=0\}. 
\end{equation}
Since ${\mathcal J}^{(i)}\subseteq\mJ^{(i+1)}$, we have 
${\mathcal J}_{(i+1)}\subseteq\mJ_{(i)}$. If $\lambda\in 
{\mathcal R}$, then $\dim \mJ^{(i)}(\lambda)=n-1+i$, which 
implies that $\dim \mJ_{(i)}(\lambda)=n-1-i$. Besides, it 
is easy to show that $ {\mathcal J}_{(i)}=V_i\oplus 
{\mathcal C}$. Therefore $\dim V_i(\lambda)=n-2-i$. In 
particular, $\dim V_{n-4}(\lambda)=2$. Also the Euler field 
$e\in V_{(n-4)}$. Fix a point $\lambda\in{\mathcal R}$. Let 
$\gamma$ be the characteristic curve, passing through 
$\lambda$, and let $\vf$ be a parametrization on $\gamma$ 
such that $\vf(\lambda)=0$. As before, let $m=n-3$. Then 
there exist a vector $\e_\vf(\lambda)\in V_{n-4}(\lambda)$
such that if $E$ and $H$ are two vector fields, satisfying 
$E\in V_{n-4}$, $H\in{\mathcal C}$, 
$E(\lambda)=\e_\vf(\lambda)$, and 
$H(\lambda)=\frac{d}{dt}\vf^{-1}(t)|_{t=0}$, 
 then
$$|\sigma\bigl(({\rm ad} H)^m E(\lambda), ({\rm ad} 
H)^{m-1} E(\lambda) \bigr)|=1.$$ Such vector is defined up 
to the transformations $\e_\vf(\lambda)\rightarrow 
\pm\e_\vf(\lambda)+\mu e(\lambda)$. 

Further, denote by $\Pi\colon\Sigma\mapsto {\mathcal R}$ 
the canonical projection. Let $\e_1$ be a vector field on 
$\Sigma$ such that 
\begin{equation}
\label{eps1} \forall u=(\lambda,\vf)\in\Sigma\quad \Pi_* 
\e_1(u)=\pm\e_\vf(\lambda)\,\,\,{\rm mod}\,\{{\mathbb R} 
e(\lambda)\}. 
\end{equation} 
Such fields $\e_1$ 
are defined modulo ${\rm span}\,\{g_0,g_1,g_2\}$ and the 
sign. How to choose among them the canonical field, up to 
the sign? 
Fix some vector field $\e_1$, satisfying (\ref{eps1}). 
Denote by $$\e_i=({\rm ad}\, h)^{i-1}\e_1,\,\,2\leq i\leq 
2m,\quad \eta=[\e_1,\e_{2m}].$$ First the tuple $(h, 
\{g_i\}_{i=0}^2,\{\e_i\}_{i=1}^{2m},\eta)$ is a frame on 
$\Sigma_D$. Let $${\mathcal L}_j={\rm span}\,\bigl\{h, 
\{g_i\}_{i=0}^2,\{\e_i\}_{i=1}^{j}\bigr\},\quad 0\leq j\leq 
2m. $$ Then one can show that 
\begin{equation*}
[\e_1,\e_2]=\kappa_1\e_2\,\, {\rm mod}\, {\mathcal L}_1 
\end{equation*}
 and, in the case $n>5$,
\begin{equation*}
[\e_1,\e_4]=\kappa_2\e_3+\kappa_3\e_4\,\,{\rm 
mod}\,{\mathcal L}_2. 
\end{equation*} 
It turns out that among all fields $\epsilon_1$, satisfying
(\ref{eps1}), there exists the unique, up to the sign,
field $\tilde\e_1$ such that the functions $\kappa_i$,
$1\leq i\leq 3$,
are identically zero. Then two frames $(h,
\{g_i\}_{i=0}^2,\{\tilde\e_i\}_{i=1}^{2m},\eta)$ and $(h,
\{g_i\}_{i=0}^2,\{-\tilde\e_i\}_{i=1}^{2m},\eta)$ are
canonically defined. This immediately implies that the
groups of symmetries is at most $(2n-1)$-dimensional.

If a $(2,n)$-distribution of maximal class has a
$(2n-1)$-dimensional group of symmetries, then all
structural functions of its canonical frames have to be
constant. It can be shown that the only nonzero commutative
relations of each of these frames in addition to the
mentioned above are
\begin{equation}\label{eq:comrel}
\begin{aligned}
~& [\tilde \e_i,\tilde \e_{2m-i+1}]=(-1)^{i+1}\eta,\,\,
[g_1, \tilde \e_i]=(2m-2i+1)\tilde \e_i,\,\, [g_2,\tilde
\e_i]=(i-1)(2m+1-i)\tilde\e_{i-1},\\ ~&[g_0,
\tilde\e_i]=-\tilde\e_i,\,\, [g_1,\eta]=2m\eta,\,\,
[g_0,\eta]=-2\eta,
\end{aligned}
\end{equation}
 which implies the uniqueness of such
distribution, up to the equivalence. Besides, from these
relations it follows that the algebra of infinitesimal
symmetries of such distribution is isomorphic to the
semi-direct sum of ${\mathfrak{gl}}(2,{\mathbb R})$
($\sim{\rm span}_{\mathbb R}\,\{ g_0,g_1,g_2,h\}$) and the
Heisenberg group ${\mathfrak n}_{2m+1}$ ($\sim {\rm
span}_{\mathbb R}\{\tilde
\e_1,\ldots,\tilde\e_{2m},\eta\}$).
Finally, it is easy to show
that for $(2,n)$-distribution, associated with
$z'(x)=y^{(n-3)}(x)$, the canonical frames satisfy the
previous commutative relations. $\Box$
\medskip

%
%
%
%
%
%

\section{Discussion}
\subsection{Distributions of non-maximal rank}
From \cite[Remark~3.4]{zelvar} it follows that a rank~2
distribution $D$ has the smallest possible class $1$ at a
point $q$ iff $\dim\, D^3(q)=4$. Suppose that $D$ satisfies
$\dim\, D^3(q)=4$ on some open set $M^o$. It is easy to
show that the distribution $D^2$ has a one-dimensional
characteristic distribution $C$. Then (locally) we can
consider the quotient $M'$ of the manifold $M^o$ by the
corresponding one-dimensional foliation together with a new
rank~2 distribution $D'$ obtained by the factorization of
$D^2$.

In fact, $D$ can be uniquely reconstructed from $D'$. Let $P(D')$ be a
submanifold in $P(TM')$ consisting of all lines lying in $D'$.
Similarly to the canonical contact system on $P(TM')$, we can define
lifts of integral curves of $D'$ to $P(D')$ and a canonical rank 2
distribution on $P(D')$ generated by tangent vectors to these lifts.
It can be proved that this contact system on $P(D')$ is locally
equivalent to $D$.

Iterating this procedure, we end up either at a non-holonomic rank 2
distribution on a three-dimensional manifold or at a distribution
$\widetilde D$, satisfying $\dim \widetilde D^3=5$. In the former case
the original distribution $D$ is locally equivalent to the Goursat
distribution and has an infinite-dimensional symmetry algebra. In
other words, the case of non-Goursat distributions of constant class~1
can be reduced to the case of distributions of class greater than~1.

This leaves the following question open: \emph{Do there exist
  completely nonholonomic rank $2$ distributions of constant class
  $2\leq m\leq n-4$?}
We know only that the answer is negative for $m=2$ ($n>5$),
which means that any such example, if it exists, should
live on at least $7$-dimensional manifold.

\subsection{Connection with Tanaka theory}

After the symplectification procedure described above, the results of this paper can be interpreted in terms of
Tanaka--Morimoto theory of structures on filtered
manifolds~\cite{mori,tan2}. The original distribution $D$
(even of maximal class) has, in general, a non-constant
symbol, which makes this theory very difficult to apply to the filtered manifold
defined by the distribution $D$ itself.
However, given rank 2 distribution $D$ of maximal class
there is a natural rank 2 distribution on the manifold
$P(\mathcal R_D)$ obtained from $\mathcal R_D$ via the
factorization by the trajectories of the Euler vector field
(or, in other words, by the projectivization of the fibers of $\mathcal
R_D$). It is generated by the projection of the sum
$V_{n-4}\oplus C$ w.r.t. this factorization. It is possible
to show that this distribution has already a fixed symbol
isomorphic to the Lie algebra generated by the vector
fields $\{h,\tilde\e_1,\dots,\tilde\e_{2m},\eta\}$ from the
proof of the main theorem (see equation~\eqref{eq:comrel}).

Moreover, there is a natural decomposition of this
distribution into the sum of two line distributions equal
to the projections of $V_{n-4}$ and $C$. This decomposition
can be interpreted as a $G$-structure on a filtered
manifold in terms of Tanaka theory and is called a
\emph{pseudo-product} structure~\cite{tan3}. The
prolongation of this structure (in terms of filtered
manifolds) is of finite type and is isomorphic to the
maximal symmetry algebra from the main theorem.

We shall dwell into the details of this approach in the forthcoming
paper.

\end{document}